\documentclass[leqno,12pt]{article}
\usepackage{amssymb,amsmath}
\usepackage{epsf,epsfig}

\newtheorem{theorem}{Theorem}
\newtheorem{proposition}[theorem]{Proposition}
\newtheorem{lemma}[theorem]{Lemma}
\newtheorem{definition}[theorem]{Definition}
\newtheorem{corollary}[theorem]{Corollary}
\newenvironment{proof}{{\bf Proof. }}{\par}

\renewcommand{\Box}
        {\hbox{\hskip 1pt \vrule width4pt height 6pt depth 1.5pt \hskip 1pt}}
\newcommand{\qed}{\relax{\ifhmode\unskip\nobreak$\kern10pt\Box$\fi
        \ifmmode\ifinner\else\hskip5pt\fi \kern10pt\Box\fi}\relax}

\newcommand{\R}{{\mathbb R}}

\renewcommand{\t}{{\mathfrak t}}

\newcommand{\CB}{{\cal B}}

\newcommand{\CH}{{\cal H}}

\newcommand{\CP}{{\cal P}}

\newcommand{\vol}{\operatorname{vol}}

\title{Arrangement of hyperplanes I: Rational functions
 and Jeffrey-Kirwan residue}
\author{Michel Brion and Mich{\`e}le Vergne}
\date{}
\begin{document}
\maketitle

\section{Introduction}

Consider the space $R_{\Delta}$ 
of rational functions of $r$ variables
with poles on an arrangement of
hyperplanes $\Delta$.
It is important to study
the decomposition of the space $R_{\Delta}$ under the action 
of the ring of differential operators with constant coefficients. 
In the one variable case, a rational function 
of $z$ with poles at most on $z=0$
is written uniquely as $\phi(z)=Princ(\phi)(z)+\psi(z)$
where $Princ(\phi)(z)=\sum_{n<0}a_n z^n$
is the principal part of $\phi(z)$ and $\psi(z)=\sum_{n\geq 0}a_n z^n$ 
is the polynomial part of $\phi(z)$.
Remark that the space 
$$G=\{\phi(z)=\sum_{n<0}a_n z^n\}$$
of principal parts 
is free under the action of $\partial/\partial z$
while the space of polynomials is evidently
a torsion module. 
Furthermore, the function $1/z$
is the unique function which cannot be written as a derivative.

We show similarly, in the case of several variables,
that there is a well determined 
decomposition of $R_\Delta$ as 
$$R_\Delta=G_{\Delta}\oplus NG_{\Delta}$$
where $G_\Delta$ is a free module under the action 
of the ring of differential operators with constant coefficients,
and $NG_\Delta$ is the torsion submodule.
Here the space $G_\Delta$  can be characterized as the space
of rational functions with a zero at infinity in all 
directions. Let us describe more precisely the space 
$G_\Delta$. We need more notations.

Let $V$ be a finite dimensional vector space
over a field $k$, of characteristic zero. 
Let $r=\dim V$.
Let $\Delta$ be a finite subset of nonzero elements of 
$V$.
Consider the  union of hyperplanes in $V^*$:
$$
{\cal H}^*(\Delta):=
\bigcup_{\alpha\in \Delta}\{z\in V^*, \langle z,\alpha\rangle=0\}
$$ 
and the ring $R_\Delta$ of rational functions on $V^*$ with poles
contained in ${\cal H}^*(\Delta)$. We denote by 
$G_{\Delta}$ the subspace of $R_\Delta$ spanned by the elements
$$
\frac{1}{\prod_{\alpha\in \kappa}\alpha^{n_\alpha}}
$$
where $\kappa$ is a subset of $\Delta$ generating $V$ and where the
$n_{\alpha}$ are positive integers.
It turns out that $G_\Delta$ is the subspace of $R_\Delta$ consisting
of functions that vanish at infinity in any direction. It is a graded
vector space with highest graded part
$G_{\Delta}[-r]:=S_{\Delta}$. Furthermore, 
$S_{\Delta}$ is  the linear span of the 
$$
\phi_\sigma=\frac{1}{\prod_{\alpha\in\sigma}\alpha}
$$
where $\sigma$ ranges over all bases of $\Delta$.

As the space $G_\Delta$ is a direct factor in $R_\Delta$, 
under the action of the ring $S(V^*)$ of differential operators with
constant coefficients, there is a natural projection $Res_\Delta$
from $R_\Delta$ to $S_\Delta$
that we call the Jeffrey-Kirwan residue.
The name Residue is justified by the fact
that the kernel of the map $Res_\Delta$ is the space of derivatives,
and by a generalization of the Cauchy formula. 
Any $S(V^*)$-morphism from $G_\Delta$
to another $S(V^*)$-module
is entirely determined by its value on $S_\Delta$,
and this morphism exists provided certain linear relations
between the $\phi_\sigma$ are satisfied. The space $S_\Delta$ is
isomorphic to the top degree component of the ``Orlik-Solomon
algebra'' associated to the hyperplane arrangement 
${\cal H}^*(\Delta)$; as a consequence, we produce bases of
$S_\Delta$ consisting of certain $\phi_\sigma$. Their dual
bases can be described in terms of iterated residues, as shown by
Szenes (see \cite{S} and section 4).

If $k=\R$, then $G_\Delta$ occurs as the space of Laplace transforms
of locally polynomial functions with possible discontinuities
on hyperplanes generated by $r-1$ elements of $\Delta$. The Laplace
transform intertwines the action of $S(V^*)$ on locally
polynomial functions by multiplication, with its action on
$G_\Delta$ by differential operators with constant coefficients.
We study the jumps of locally polynomial functions
in terms of the poles of their Laplace transforms. As a
consequence, we show that a locally polynomial function is continuous
if and only if its Laplace transform vanishes at order 2 in any
direction. We also construct inverses of the Laplace transform, using
our description of $G_\Delta$ by generators and relations.

Many of the statements proved in this article are already implicitly
stated in Jeffrey-Kirwan articles \cite{JK1} and \cite{JK2}.
However, we felt the need, for applications, to clarify
some of their statements. The main application will be 
an algebraic construction of Eisenstein series: to each rational
function with poles on hyperplanes, we will associate 
a periodic meromorphic function in several variables. This will be
treated in part II of this article.

Applications to the Poisson formula will be given in another article.

Our interest in the space of functions $R_\Delta$ and
their Laplace transforms comes from the study of 
integrals over symplectic spaces of equivariant cohomology classes.
Let $(M,\Omega)$ be a compact symplectic
manifold, with an Hamiltonian action of a torus $T$. 
Let $f:M\to\t^*$ be the moment map. Let $X\in \t$.
Let $\Omega(X)=\langle f,X\rangle+\Omega$ be the equivariant
symplectic form. Let $\alpha(X)$ be an equivariant closed form on $M$.
 Consider the integral 
$$I(X)=\int_M \alpha(X) e^{\Omega(X)}.$$
Assume for simplicity that the set $F$ of fixed points for 
the action of $T$ on $M$ is finite. For $p\in F$, 
let $\Delta_p\subset\t^*$ be the set of weights 
for the action of $T$ in the tangent space $T_pM$.
Then, by the localisation formula in equivariant cohomology,
we have 
$$I(X)=\sum_{p\in F}\phi_p(X)e ^{\langle f(p),X\rangle}$$
where each $\phi_p$ is in the ring $R_{\Delta_p}.$
 
If $\xi$ is a regular value of $f$, we can consider 
the reduced space $M_{red}(\xi)=f^{-1}(\xi)/T$
with reduced symplectic structure $\Omega_\xi$.
The equivariant cohomology class of $\alpha(X)$
gives rise to a de Rham cohomology class $\alpha_\xi$ on
$M_{red}(\xi)$. Consider the function 
$$r(\xi)=\int_{M_{red}(\xi)}\alpha_\xi e ^{\Omega_\xi}.$$
This function is defined for regular values of $\xi$.
It is important to determine this function
and its jumps when crossing walls of singular values of the moment
map. The functions $I(X)$ and $r(\xi)$ are related by the Laplace
transform. Thus it is important to study jumps of Laplace transforms
of functions in the space $G_\Delta$.

We thank Michel Duflo for his comments on this article,
and Bernard Malgrange for his decisive help in the proof 
of Theorem 1.


\section{Rational functions with poles on hyperplanes:
Jeffrey-Kirwan residue}

Let $V$ be a finite dimensional vector space
over a field $k$, of characteristic zero. 
Let $r=\dim V$.
We denote by $S(V)$ the symmetric algebra of $V$.
Let $V^*$ be the dual space.
We identify $S(V)$ with the ring of polynomial functions 
on $V^*$.
Let $\Delta\subset V$ be a finite subset of nonzero elements,
which spans $V$.
We denote by 
$$R_\Delta:=\Delta^{-1}S(V)$$ 
the ring generated over $S(V)$ by inverting the linear functions 
$\alpha\in \Delta$. This is a ring graded by the degree (positive 
or negative). Consider the  union of hyperplanes in $V^*$:

$$
\CH^*(\Delta):=
\bigcup_{\alpha\in \Delta}\{z\in V^*, \langle\alpha,z\rangle=0\}
$$ 
and the  open subset 
$$V^*_{reg,\Delta}=V^*-\CH^*(\Delta)$$
of ($\Delta$)-regular elements in $V^*$. Then $R_\Delta$
is the ring of rational functions on $V^*$ with poles contained 
in the union  of hyperplanes $\CH^*(\Delta)$. Functions in 
$R_\Delta$ are defined on the set $V^*_{reg,\Delta}$ of
regular elements.

Let ${\cal D}$ be the ring of differential operators on $V^*$ with
polynomial coefficients. Recall that the ring ${\cal D}$ is
generated by its subrings $S(V)$ of polynomial functions on $V^*$, and
$S(V^*)$ of differential operators on $V^*$ with constant
coefficients. Observe that $S(V)$ and $R_\Delta$ are 
graded ${\cal D}$-modules.

If $\phi\in R_{\Delta}$ and if $y\in V^*$ is a regular element,
then $t\mapsto \phi(y+tz)$ is a rational function for any 
$z\in V^*$. We say that $\phi$ {\bf vanishes at infinity}
if the rational function $t\mapsto \phi(y+tz)$ is $0$ at 
$\infty$ for all regular $y\in V^*$ and for all $z\in V^*$.

Let $\kappa$ be a subset of $\Delta$.
The subset $\kappa$ is called {\bf generating} if the 
$\alpha\in \kappa$  generate the vector space $V$.
It is called a {\bf basis} of $\Delta$, if the $\alpha\in \kappa$ 
form a basis of $V$.
We denote by ${\cal B}(\Delta)$
the set of bases  of $\Delta$.
 
For $\kappa\subset \Delta$, set 
$$
\phi_\kappa:=\frac{1}{\prod_{\alpha\in \kappa}\alpha}.
$$

We denote by $G_{\Delta}$ the subspace of 
$R_\Delta$ spanned by the 
$$\frac{1}{\prod_{\alpha\in \kappa}\alpha^{n_\alpha}}$$
where $\kappa$ is {\bf generating}
and the $n_\alpha$ are  positive integers.
Then $G_{\Delta}$ is a graded vector space with highest graded part 
$S_{\Delta}$ (in degree $-r$). Furthermore,
$S_{\Delta}$ is  the linear span of the 
$\phi_\sigma$ where $\sigma$ ranges over all bases of 
$\Delta$.

Clearly, any function in the space $G_\Delta$ vanishes at infinity. We
will prove that the converse holds in Theorem 1 below.
\bigskip

\noindent
{\bf Remark.} The space $G_\Delta$ is contained in 
$\sum_{j\geq r}R_{\Delta}[-j]$
but is strictly smaller if $r>1$. 
For example, if $\alpha\in \Delta$, then 
$\alpha^{-r}$ is never in $G_\Delta$.
\bigskip

We denote by 
$NG_{\Delta}$ the subspace 
of $R_{\Delta}$ spanned by the 
$$
\frac{\psi}{\prod_{\alpha\in \kappa}\alpha^{n_\alpha}}
$$
where $\psi\in S(V)$, $\kappa$ is {\bf not generating}
and the $n_\alpha$ are non-negative integers.

Remark that the subspace $NG_{\Delta}$ of $R_{\Delta}$ is stable 
under the action of ${\cal D}$, whereas $G_{\Delta}$ is stable under
the action of $S(V^*)$ by differential operators with constant
coefficients.
\begin{theorem}\label{free}
We have a direct sum decomposition
of $S(V^*)$-modules
$$
R_{\Delta}=G_{\Delta}\oplus NG_{\Delta}.
$$
Moreover, the space $G_{\Delta}$ is a free 
$S(V^*)$-module, and is freely  generated by $S_{\Delta}$, 
while the space $NG_{\Delta}$ is the torsion submodule. Finally,
$G_{\Delta}$ is the space of functions in $R_{\Delta}$ 
which vanish at infinity.
\end{theorem}

For this we prove a succession of lemmas.
\begin{lemma}\label{lemma1}
The $S(V^*)$-module $G_{\Delta}$ is generated by 
$S_{\Delta}$.
Moreover, we have $R_{\Delta}=G_{\Delta}+NG_{\Delta}$.
\end{lemma}
\begin{proof}
Observe that the $S(V^*)$-module generated by $S_{\Delta}$ is the
span of the elements
$$\frac{1}{\prod_{\beta\in\sigma} \beta^{n_\beta}}$$
where $\sigma\in{\cal B}(\Delta)$ and where each $n_{\beta}$ is a
positive integer. To prove the first assertion, it is enough to check
that this vector space is stable by multiplication by
$1/\alpha^n$ where $\alpha\in\Delta$. For this, write
$\alpha=\sum_{\beta\in\sigma} c_{\alpha\beta}\beta$.
Then we have
$$
\frac{1}{\alpha^n\prod_{\beta\in\sigma}\beta^{n_{\beta}}}=
\sum_{\beta\in\sigma} \frac{c_{\alpha\beta}}
{\alpha^{n+1}\beta^{n_{\beta}-1}
\prod_{\gamma\in\sigma,\gamma\neq \beta}\gamma^{n_{\gamma}}}.
$$
If $\beta\in\sigma$ is such that $n_{\beta}=1$, then the 
corresponding term in the right-hand side is in the 
$S(V^*)$-module generated by $G_{\Delta}[-r)$: indeed, 
if $c_{\alpha\beta}\neq 0$ then
$\sigma\cup\{\alpha\}\setminus\{\beta\}$ is a basis of $\Delta$.
On the other hand, if $n_{\beta}>1$ then our term is the inverse of
$\alpha^{n+1}\prod_{\beta\in\sigma} \beta^{n'_{\beta}}$ with
$n'_{\beta}\geq 1$ and 
$\sum_{\beta\in\sigma}n'_{\beta}=(\sum_{\beta\in\sigma}n_{\beta})-1$.
So the assertion follows by induction on 
$\sum_{\beta\in\sigma} n_{\beta}$.

Similarly, any element of $R_\Delta=\Delta^{-1}S(V)$ is a linear
combination of elements
$$
\phi=\frac{\psi}{\prod_{\alpha\in\kappa}\alpha^{n_{\alpha}}}
$$
where $\psi\in S(V)$, $\kappa$ is linearly independent and the
$n_{\alpha}$ are positive integers. If moreover $\kappa$ is not
generating, then $\phi$ is in $NG_{\Delta}$. If $\kappa$ is generating,
then we can express $\psi$ as a polynomial in the variables
$\alpha\in\kappa$, and we obtain $\phi\in G_{\Delta}+NG_{\Delta}$.
\end{proof}
\bigskip

\begin{lemma} The $S(V^*)$-module $R_{\Delta}/NG_{\Delta}$ is free.
\end{lemma}
\begin{proof}
Observe that $R_{\Delta}/NG_{\Delta}$ is a ${\cal D}$-module.
Furthermore, it is spanned (as a vector space) by the images of 
$$
\frac{1}{\prod_{\alpha\in\sigma}\alpha^{n_{\alpha}}}
$$
where $\sigma$ is a basis of $\Delta$, and the $n_{\alpha}$ are
positive integers. It follows that the ${\cal D}$-module
$R_{\Delta}/NG_{\Delta}$ is generated by the images
$\overline{\phi_{\sigma}}$ of the $\phi_{\sigma}$ 
($\sigma\in{\cal B}(\Delta)$). Observe that
$\overline{\phi_{\sigma}}$ is killed by $V$; thus, the 
${\cal D}$-module ${\cal D}\overline{\phi_{\sigma}}$ is a non zero 
quotient of ${\cal D}/{\cal D}V$. The latter is a simple 
${\cal D}$-module, isomorphic to $S(V^*)$; therefore, 
${\cal D}\overline{\phi_{\sigma}}$ is isomorphic to $S(V^*)$, too.
Iterating this argument, we construct an ascending filtration of 
the ${\cal D}$-module $R_{\Delta}/NG_{\Delta}$, each submodule being 
generated by certain $\overline{\phi_{\sigma}}$'s, with successive 
quotients isomorphic to $S(V^*)$.
\end{proof}

\begin{lemma} The subspace 
$S_{\Delta}$ intersects $NG_{\Delta}$ trivially.
\end{lemma}
\begin{proof}
We argue by induction on the number of elements in
$\Delta$. We may assume that 
$\Delta$ contains no proportional elements. Let 
$\phi\in S_{\Delta}$.
Write 
$$
\phi=\sum_{\sigma\in{\cal B}(\Delta)} 
\frac{c_{\sigma}}{\prod_{\alpha\in\sigma} \alpha}
$$
and consider $\phi$ as a rational function on $V^*$. 
Observe that the poles of $\phi$ are simple and 
along the hyperplanes $\alpha=0$ ($\alpha\in\Delta$).
Choose $\alpha$ among the poles of $\phi$. Choose a
decomposition $V=k\alpha\oplus V_0$. Then $Q(V)$ (the fraction field
of $S(V)$) is identified with the field of rational functions in the
variable $\alpha$, with coefficients in $Q(V_0)$. Therefore, we have a
restriction map $S(V)\to S(V_0):\phi\mapsto\phi_0$.
Consider the image $\Delta_0$ of $\Delta\setminus\{\alpha\}$ in $V_0$.
The restriction map extends to an homomorphism
$(\Delta\setminus \{\alpha\})^{-1}
S(V)\to \Delta_0^{-1}S(V_0)$ by restriction to generic points.
We have also a residue map 
$Res_{\alpha}:Q(V)\to Q(V_0)$ 
with respect to the variable $\alpha$, defined by the formula 
$$
Res_{\alpha}(\phi)= \frac{1}{(K-1)!}
\left((\frac{\partial}{\partial\alpha})^{K-1}(\alpha^{K}\phi)\right)_0
$$
for any integer $K$ such that 
$\alpha^K\phi\in R_{\Delta\setminus\{\alpha\}}$.
 
As $\alpha$ is a simple pole of $\phi$, we have simply 
$$
Res_{\alpha}(\phi)= \sum_{\sigma,\alpha\in\sigma} 
\frac{c_{\sigma}}{\prod_{\beta\in\sigma,\beta\neq\alpha}\beta_0}
$$
where $\beta_0$ denotes the image of $\beta$ in $V_0$.
If $\sigma$ is a basis of $\Delta$ which contains $\alpha$, then
$(\sigma\setminus\{\alpha\})_0$ is a basis of $\Delta_0$. Therefore,
$Res_{\alpha}(\phi)$ is in $G_{\Delta_0}$.

Consider a generator
$$
u=\frac{\psi}{\prod_{\beta\in \kappa}\beta^{n_\beta}}
$$
of $NG_\Delta$, with $\psi\in S(V)$ and $\kappa$ non generating.
Write 
$$
u= \frac{\psi}
{\alpha^K\prod_{\beta\in \kappa, \beta\neq \alpha}\beta^{n_\beta}}.
$$
If $K=0$, then $Res_{\alpha}(u)=0$. 
If $K>0$, the set $\kappa$ contains $\alpha$ and 
is non generating. Thus, its restriction $\kappa_0$ is non generating.
We see that $Res_{\alpha}(u)$ can be written as  
$$
Res_{\alpha}(u)=
\frac{\psi'}{\prod_{\beta_0\in\kappa_0} \beta_0^{n_\beta+K-1}}.
$$
for some $\psi'\in S(V_0)$,
so that $Res_{\alpha}(u)\in NG_{\Delta_0}$.

If $\phi\in G_{\Delta}\cap NG_{\Delta}$,
it follows from the above discussion that 
$Res_{\alpha}(\phi)\in G_{\Delta_0}\cap NG_{\Delta_0}$.
Therefore, by the induction hypothesis, we have
$Res_{\alpha}(\phi)=0$: thus, $\phi$ has no pole along $\alpha=0$. By
the beginning of the proof, $\phi$ has no pole at all, so that
$\phi=0$. 
\end{proof}
\begin{lemma}
If $\phi\in R_{\Delta}$ vanishes at infinity, then $\phi$ is in
$G_{\Delta}$.
\end{lemma}
\begin{proof}
First we claim that the space of functions which vanish at infinity is
stable by the action of $S(V^*)$. Indeed, let $\phi\in R_{\Delta}$
vanish at infinity. Write
$$
\phi=\frac{\psi}{\prod_{\alpha\in\Delta}\alpha^{n_{\alpha}}}
$$
where $\psi\in S(V)$. For $z\in V^*$, set
$$
n(z):=\sum_{\alpha,\langle\alpha,z\rangle\neq 0}\,n_{\alpha}.
$$
The assumption that $\phi$ vanishes at infinity means that
$$
deg(t\mapsto \psi(y+tz)) < n(z)
$$
for all regular $y$ and for all $z$ in $V^*$. Let $w\in V^*$; then, for
all $u\in k$ such that $y+uw$ is regular, we also have
$$
deg(t \mapsto \psi(y+tz+uw)) < n(z)
$$
and therefore, the function
$$
\frac{\partial(w)\psi}{\prod_{\alpha\in\Delta}\alpha^{n_{\alpha}}}
$$
vanishes at infinity. Now
$$\partial(w)\phi=
\frac{\partial(w)\psi}{\prod_{\alpha\in\Delta}\alpha^{n_{\alpha}}}
-\sum_{\alpha\in\kappa}
\frac{n_{\alpha}\langle\alpha,w\rangle}{\Delta}\phi
$$
which implies our claim.

Assume now that there exists a non-zero $\phi\in NG_{\Delta}$ which
vanishes at infinity. As in the proof of Lemma \ref{lemma1}, we can
write
$$
\phi=\sum_{\kappa}\,\phi_{\kappa}
$$
where the sum is over all linearly independent subsets
$\kappa\subset\Delta$ which are not bases, and where each
$\phi_{\kappa}$ is in $\kappa^{-1}S(V)$. Furthermore, we may assume
that the number of $\kappa$ such that $\phi_{\kappa}$ is non-zero is
minimal (among all possible decompositions of all non-zero 
$\phi\in NG_{\Delta}$ which vanish at infinity).

Choose $\kappa_0$ such that $\phi_{\kappa_0}\neq 0$, and choose
a non-zero $z_0\in V^*$ such that $\langle\alpha,z_0\rangle=0$ for all
$\alpha\in\kappa_0$. Then $\partial^n(z_0)\phi_{\kappa_0}=0$
for large $n$. But all successive derivatives of $\phi$ vanish at
infinity and are in $NG_{\Delta}$. Moreover,
$$
\partial^n(z_0)\phi=\
\sum_{\kappa\neq\kappa_0}\,\partial^n(z_0)\phi_{\kappa}
$$
is a decomposition with fewer terms than $\phi$. Thus,
$\partial^n(z_0)\phi=0$ for some positive $n$.

Choose $n$ minimal with this property, and set 
$\psi:=\partial^{n-1}(z_0)\phi$. Then $\psi$ is a non-zero element of
$NG_{\Delta}$ which vanishes at infinity, and $\partial(z_0)\psi=0$.
But then the function $t\mapsto\psi(y+tz_0)$ is constant for any 
$y\in V^*$, a contradiction.
\end{proof}
\bigskip

The space $S_\Delta$
is generated by the elements $\phi_\sigma$ where 
$\sigma$ ranges over $\CB(\Delta)$. However, there 
are linear relations among the elements $\phi_\sigma$.
Indeed, let $\sigma$ be a basis
of $\Delta$. If 
$\alpha\in\Delta\setminus\sigma$ and
$$
\alpha=\sum_{\beta\in\sigma} c_{\alpha\beta}\beta
$$
is the expansion of $\alpha$ in the basis $\sigma$,
then $\sigma\cup\{\alpha\}\setminus\{\beta\}$
is a basis if and only if $c_{\alpha\beta}$
is non zero, and we have  
$$
\phi_{\sigma}=\sum_{\beta\in\sigma, c_{\alpha\beta}\neq 0} c_{\alpha\beta}
\phi_{\sigma\cup\{\alpha\}\setminus\{\beta\}}.
$$ 
In section 4, we will prove that the linear relations between the
elements $\phi_\sigma$ are generated by the relations above.

\bigskip

We can now define the Jeffrey-Kirwan residue map: denote by 
$$
\hat R_{\Delta}:=\Delta^{-1}\hat S(V)
$$
the ring  generated over the ring $\hat S(V)$ 
of formal power series, by inverting the linear 
functions $\alpha\in \Delta$.
Define the Taylor expansion at order $K$
as the projection  
$$
Taylor_{[\leq K]}: \hat R_{\Delta}\to\bigoplus_{j\leq K}R_{\Delta}[j].
$$
Using $Taylor_{[\leq -r]}$, we project the space
$\hat R_{\Delta}$ to $R_{\Delta}[\leq -r]$.
Then using the direct sum decomposition
$$
R_{\Delta}=G_{\Delta}\oplus NG_{\Delta}
$$
we obtain a projection map
$$
Princ_{\Delta}: \hat R_{\Delta}\to G_{\Delta}
$$
by composing both projections 
$\hat R_{\Delta}\to R_ {\Delta}[\leq -r]\to  G_{\Delta}$.

Remark that as $G_{\Delta}$ is contained in $R_{\Delta}[\leq -r]$,
the map $Princ_{\Delta}$
can also be defined as the composition of  
$Taylor_{[\leq K]}:\hat R_{\Delta}\to R_ {\Delta}[\leq K]$
for any index $K\geq -r$,
followed by the  projection $R_{\Delta}[\leq -K]\to G_{\Delta}$.
\begin{definition}
The Jeffrey-Kirwan residue map
$$Res_{\Delta}:\hat R_{\Delta}\to S_{\Delta}$$ 
is defined to be the composite of the projection
$Princ_{\Delta}$ followed by the projection of $G_{\Delta}$ on 
$S_{\Delta}$.
\end{definition}

In other words, the map 
$Res_{\Delta}$ is the identity on 
$S_{\Delta}$ and vanishes on 
$\oplus_{j\neq -r}R_{\Delta}[j]$ and on $NG_{\Delta}$ as well.
We can determine easily the map 
$Res_{\Delta}$ on 
$\hat R_{\Delta}$ 
by first projecting
on $R_{\Delta}[-r]$, then, using the fact that 
$Res_{\Delta}$ vanishes on 
$NG_{\Delta}$, projecting further on $G_{\Delta}[-r]=S_{\Delta}$.

Consider the subspace $V^*R_\Delta$ spanned by derivatives of elements
of $R_\Delta$; it is a submodule of $R_{\Delta}$ under the action of
$S(V^*)$.

\begin{proposition}\label{ker}
We have
$$
V^*R_{\Delta}=NG_{\Delta}\oplus\bigoplus_{j<-r}G_{\Delta}[j].
$$
In particular, we have
$$
R_{\Delta}=V^*R_{\Delta}\oplus S_{\Delta},~
\hat R_{\Delta}=V^*\hat R_{\Delta}\oplus S_{\Delta}
$$
and the kernel of $Res_{\Delta}$ is $V^*\hat R_{\Delta}$.
\end{proposition}
\begin{proof}
From Theorem \ref{free} we obtain
$$
V^*R_{\Delta}=V^*NG_{\Delta} \oplus V^*G_{\Delta}
=V^*NG_{\Delta} \oplus \bigoplus_{j<-r} G_{\Delta}[j].
$$
So it is sufficient to check that 
$NG_{\Delta}=V^*NG_{\Delta}$. For this, consider
$$\phi=\frac{\psi}{\prod_{\alpha\in\kappa} \alpha^{n_{\alpha}}}$$
where $\psi\in S(V)$ and where $\kappa$ is linearly dependent. Choose
$y\in V^*$ such that $\langle y,\alpha\rangle=0$ for all
$\alpha\in\kappa$. We can find $\Psi\in S(V)$ such that
$\partial(y)\Psi=\psi$; then 
$$
\phi=\partial(y)(\frac{\Psi}{\prod_{\alpha\in\kappa}\alpha^{n_{\alpha}}}).
$$
\end{proof}
\bigskip

In particular, the kernel of $Res_{\Delta}$ is the space of
derivatives. Using $Res_{\Delta}$, we now obtain a multidimensional
analogue of the Cauchy formula: for any meromorphic function $\phi$ of
one variable $z$, and for any $y\neq 0$, we have 
$$(Princ\,\phi)(y)=Res_{z=0}\frac{\phi(z)}{y-z}.$$

Let $y\in V^*$ be regular and let $\psi\in R_{\Delta}$. Set
$$
(C(y)\psi)(z):=\psi(y-z).
$$
Then the rational function $C(y)\psi$ is defined at $0$, and thus its
Taylor series at the origin is in $\hat S(V)$. 
To any $\phi\in\hat R_{\Delta}$, we associate the
endomorphism $u(\phi)$ of $S_{\Delta}$ defined by 
$$(u(\phi)\psi)(z) = Res_{\Delta}(\phi(z)\psi(y-z)).$$
Denoting by $m(\phi)$ the multiplication by $\phi$, then 
$u(\phi)$ is composition $Res_{\Delta}\circ m(\phi)\circ C(y)$. 
We consider its trace 
$Tr_{S_\Delta} (Res_\Delta \circ m(\phi)\circ C(y))$.
\begin{proposition}\label{cauchy}
For any regular $y$ in $V^*$ and for any $\phi\in\hat R_{\Delta}$, we
have 
$$(Princ_{\Delta}\phi)(y)=
Tr_{S_{\Delta}} (Res_{\Delta}\circ m(\phi)\circ C(y)).$$
\end{proposition}
\begin{proof}
First we consider the case where $\phi\in S_{\Delta}$. Then 
$$Res_{\Delta}(\phi(z)\psi(y-z))=\phi(z)\psi(y)$$
because $z\mapsto\psi(y-z)$ is defined at $0$. So $u(\phi)$
maps $\psi$ to $\psi(y)\phi$, and its trace is 
$\phi(y)=Princ_{\Delta}(\phi)(y)$.

Now we assume that the formula holds for $\phi$, and we claim that it
holds for $\partial(w)\phi$ where $w\in V^*$. Indeed, using the fact
that $Res_{\Delta}$ vanishes on derivatives, we obtain
$$\displaylines{
Res_{\Delta}((\partial(w)\phi)(z)\psi(y-z))
=-Res_{\Delta}(\phi(z)\partial_z(w)\psi(y-z))
\hfill\cr\hfill
=Res_{\Delta}(\phi(z)\partial_y(w)\psi(y-z))
=\partial_y(w)Res_{\Delta}(\phi(z)\psi(y-z))
\hfill\cr}$$
which implies the claim.

It follows that the formula holds for any $\phi\in G_{\Delta}$. If
$\phi\in\widehat{NG_{\Delta}}$ then the left-hand side vanishes. On
the other hand, the function $z\mapsto\psi(y-z)$  is in 
$\hat S(V)$; thus, $z\mapsto \phi(z)\psi(y-z)$ is in 
$\widehat{NG_{\Delta}}$ and the right-hand side vanishes, too.
\end{proof}
\bigskip

\noindent
{\bf Remark.} More generally, let 
$A:\hat R_{\Delta}\to\hat R_{\Delta}$ be an operator which commutes
with the action of $S(V^*)$. Then we have for any regular $y\in V^*$
and for any $\phi\in \hat R_{\Delta}$:
$$
A(Princ_{\Delta}(\phi))(y)=
Tr_{S_{\Delta}}(Res_{\Delta}\circ m(\phi)\circ C(y)\circ A)
$$
(the proof is the same).
\bigskip

Let us deduce from this (abstract) Cauchy formula, an explicit
expression of $Princ_{\Delta}(\phi)$ in terms of derivatives of
elements of $S_{\Delta}$. For this, choose a basis 
$(\phi_b)_{b\in B}$ of $S_{\Delta}$ and denote by $(\phi^b)$ the
dual basis. For $\phi\in\hat R_{\Delta}$ and $h\in V$, the function 
$y\mapsto e^{-\langle y,h\rangle}\phi(y)$ is in $\hat R_{\Delta}$.
Moreover, the map
$$
h\mapsto \langle\phi^b,Res_{\Delta}(e^{-h}\phi)\rangle
:=D^b(\phi)(h)
$$
is easily seen to be polynomial. It thus defines a differential
operator $D^b(\phi)$ on $V^*$.
\begin{proposition}
For any $\phi\in \hat R_{\Delta}$,
and for any basis $(\phi_b)_{b\in B}$ of $S_{\Delta}$,
we have
$$Princ_{\Delta}(\phi)=\sum_{b\in B} D^b(\phi)\cdot \phi_b.$$
\end{proposition}
\begin{proof}
Let $y$ be a regular element of $V^*$. Then we have by the Cauchy
formula:
$$\displaylines{
Princ_{\Delta}(\phi)(y)=
Tr_{S_{\Delta}}(Res_{\Delta}\circ m(\phi)\circ C(y))
\hfill\cr\hfill
=\sum_{b\in B}
\langle\phi^b,Res_{\Delta}(\phi(z)\phi_b(y-z))\rangle.
\cr}$$
Now observe that $\phi_b(y-z)=(e^{-\partial(z)}\phi_b)(y)$. Thus, we
have
$$\langle \phi^b,Res_{\Delta}(\phi(z)\phi_b(y-z))\rangle
=D^b(\phi)\cdot\phi_b.$$
\end{proof}

\bigskip
Remark that Propositions \ref{basis} and \ref{iterated} below
provide a basis $(\phi_b)_{b\in B}$ together with the dual basis 
$(\phi^b)_{b\in B}$.
Thus we obtain an explicit expression of any element 
in $G_{\Delta}$ as a sum of successive derivatives of elements
$\phi_{\sigma}$. This provides a way of separating variables.

\noindent
{\bf Example.} Let $V$ be a vector space with basis $(e_1,e_2)$.
Let $\Delta$ be the ordered set
$$\Delta=(e_1, e_2, e_1+e_2).$$
The set $B$ of Proposition \ref{basis}
according to this ordering consists of 
$$b_1=(e_1,e_2)\hspace{1cm} b_2=(e_1, e_1+e_2).$$
Furthermore, if $\sigma=\{e_2, e_1+e_2\}$,
we have $\phi_\sigma=\phi_{b_1}-\phi_{b_2}$. Let 
$$
\phi(z_1,z_2)=\frac{1}{z_1 z_2 (z_1+z_2)}.
$$
If $h=h_1 e_1+h_2 e_2$, 
the component of degree $-2$ of 
$e^{-h_1z_1-h_2z_2}\phi(z_1,z_2)$
is $$\frac{-h_1z_1-h_2z_2}{z_1z_2(z_1+z_2)}
=-\frac{h_1}{z_2 (z_1+z_2)}-\frac{h_2}{z_1 (z_1+z_2)}
=-h_1 \phi_{b_1}+(h_1-h_2) \phi_{b_2}.$$

We have indeed  
$$\frac{1}{z_1 z_2 (z_1+z_2)}=
-\frac{\partial}{\partial z_1}
\cdot \frac{1}{z_1 z_2}+
(\frac{\partial}{\partial z_1}-\frac{\partial}{\partial z_2})
\cdot \frac{1}{z_1(z_1+z_2)}.$$
\bigskip

\noindent
{\bf Remark.} The residue that Jeffrey and Kirwan actually defined
is a linear form over $S_\Delta$, defined in the case  when
$k=\R$. It depends on choices of chambers in $V$ and $V^*$. We will
describe this residue in section 5.


\section{Residue along a hyperplane} 

Let us recall the notion of a residue map along a hyperplane.

Let $V_0$ be an hyperplane in $V$.  We denote
by $\Delta_0$ the subset $\Delta\cap V_0$.
The space  $V_0^{\perp}$ is a line in $V^*$.
The fibers of the restriction map 
$V^*\to V_0^*$ are affine lines $z+V_0^{\perp}$.
If $\phi$ is a rational function with poles on the set of hyperplanes
$\Delta$,
its restriction to the affine line 
$z+V_0^{\perp}$
is a rational function, except when the affine line $z+ V_0^{\perp}$
is contained in the pole set of $\phi$
(in this case the restriction is nowhere defined).
The residue at infinity of this rational function
is well defined.
More precisely, choose diffential forms of maximal degree
$\omega$ on $V^*$, $\omega_0$ on $V_0^*$
and choose an equation $z_0$ of $V_0$, 
such that $\omega_0={\rm int}(z_0)\omega$ where ${\rm int}$ is the
contraction. Define the residue map
$$Res_{V/V_0}:\Delta^{-1}S(V)\otimes\wedge^r V\to
\Delta_0^{-1}S(V_0)\otimes\wedge^{r-1}V_0
$$
by
$$
Res_{V/V_0}(\phi\otimes \omega)(z)=-Res_{t=\infty}
(\phi(z+tz_0)dt)\otimes \omega_0
$$
for $z\in V^*$ (clearly, this only depends on the image of $z$ in
$V^*/kz_0=V_0^*$). We now give a characterization of this map. 

We identify $R_{\Delta_0}=\Delta_0^{-1}S(V_0)$ to a subalgebra of 
$R_{\Delta}$, so that
$R_{\Delta}$ is a $R_{\Delta_0}$-module. We denote by 
$\Delta_1$ the complement of $\Delta_0$ in $\Delta$.

If $\nu=(\alpha_j, 1\leq j\leq L)$
is a sequence of elements of $\Delta$
with possible repetitions, we set
$$
m_\nu:=\frac{1}{\prod_{j=1}^L\alpha_j}.
$$
We write $\nu\subset\Delta_0$ (resp. $\nu\subset\Delta_1$)
if all elements $\alpha_j$ of the sequence $\nu$ are in $\Delta_0$
(resp. $\Delta_1$).

Let $\omega\in \wedge^rV$.
If $\beta\in \Delta_1$, there exists a unique 
$\omega_0\in \wedge^{r-1}V_0$ such that 
$\omega= \beta\wedge \omega_0$.
We then write 
$$\omega_0= (d\beta)^{-1}\wedge \omega.
$$
\begin{proposition}\label{existence}
The map $Res_{V/V_0}$
is the  unique
$R_{\Delta_0}$-linear map 
$$
res_{V/V_0}:R_{\Delta}\otimes \wedge^rV
\to R_{\Delta_0}\otimes \wedge^{r-1}V_0
$$
such that, for $\omega\in \wedge^rV$,

1) for any $\beta\in \Delta_1$,
$$
res_{V/V_0}( \frac{1}{\beta}\otimes\omega)
=(d\beta)^{-1}\wedge \omega.
$$

2) 
$$
res_{V/V_0}(S(V)\otimes\omega)=0.
$$

3) for any sequence $\nu\subset \Delta_1$:
$$
res_{V/V_0}(m_{\nu}\otimes\omega)=0
$$
if the length of $\nu$ is strictly greater than $1$.
\end{proposition}
 
Indeed, these properties are easily checked for the map $Res_{V/V_0}$
defined above, and uniqueness follows from the following remark.
\begin{proposition}\label{true}
We have 
$$R_{\Delta}=\Delta_0^{-1}S(V)+\sum_{\nu\subset \Delta_1}
R_{\Delta_0}m_{\nu}.$$
\end{proposition}
\begin{proof}
Let $\psi\in S(V)$ and $\nu$ a sequence of elements of $\Delta$.
Consider the element $\psi m_{\nu}$ of $R_{\Delta}=\Delta^{-1}S(V)$. 
If $\nu$ is contained in $\Delta_0$, or if $\psi\in S(V_0)$, we are
already in the desired set. If $\alpha_j\in\nu$ is not in $\Delta_0$
and if $\psi$ is not in $S(V_0)$, then using the decomposition 
$$
S(V)=S(V_0)\oplus\alpha_j S(V)
$$
we can strictly decrease the power of 
$\alpha_j$ in the expression of $m_\nu$.
\end{proof} 
\bigskip

We finally note some properties of the map $Res_{V/V_0}$.

We extend the map $Princ_{\Delta}:R_\Delta\to G_\Delta$
to a map 
$$
Princ_{\Delta}:R_\Delta\otimes \wedge^{max}V\to
G_\Delta\otimes \wedge^{max}V
$$ 
still denoted by $Princ_{\Delta}$.
In the same way we extend the map $Res_{\Delta}$
to a map 
$$
Res_{\Delta}:R_\Delta\otimes \wedge^{max}V
\to S_\Delta\otimes \wedge^{max}V.
$$
\begin{proposition}
The map 
$Res_{V/V_0}$ is homogeneous of degree $-1$, and
is compatible with the maps $Princ$
and with the Jeffrey-Kirwan residue.
More explicitly:

1) 
$$
Res_{V/V_0} (Princ_{\Delta}(\phi))=
Princ_{\Delta_0} (Res_{V/V_0} (\phi)).
$$

2) 
$$
Res_{V/V_0} (Res_{\Delta}(\phi))=
Res _{\Delta_0} (Res_{V/V_0}(\phi)).
$$
\end{proposition}
\begin{proof}
Remark that  
$$
NG_{\Delta}\subset\Delta_0^{-1}S(V) +
\sum_{\nu\subset \Delta_1}\,NG_{\Delta_0}m_{\nu}.
$$
Therefore, $Res_{V/V_0}$ maps $NG_{\Delta}$ to $NG_{\Delta_0}$, 
and both members of equation $(1)$ vanish on $NG_\Delta$.

Now consider an element $m_\nu=m_{\nu_0}m_{\nu_1}$
where $\nu$ is generating.
If the length of $\nu_1$ is greater than $1$, both members of equation 
(1) vanish. If $\nu_1$ consists of one element, then $\nu_0$
generates $V_0$ and we obtain Assertion 1. Assertion 2 follows from
the fact that $Res_{V/V_0}$ is homogeneous of degree -1.
\end{proof}

\section{Orlik-Solomon relations}

In this section, we describe the linear relations between the
generators $\phi_\sigma$ ($\sigma\in{\cal B}(\Delta)$) of the space
$S_\Delta$, and we construct bases of this space consisting of
certain $\phi_\sigma$. Using iterated residues, we construct the
dual bases as well.

For this, we begin by interpreting the space $S_\Delta$ in terms
of the {\bf Orlik-Solomon algebra} associated to the hyperplane
arrangement ${\cal H}^*(\Delta)$, see \cite{OT} Chapter 3. Recall that
this algebra, which we denote by $A_{\Delta}$, is the subalgebra of
rational differential forms on $V^*$ generated by the forms
$$
\omega_\alpha:=\frac{d\alpha}{\alpha}
$$
where $\alpha\in\Delta$. Clearly, $A_\Delta$ is graded by the degree
of differential forms, and its top degree component is
$$
A_\Delta[r]=S_\Delta\otimes\wedge^r V.
$$

It is known that the algebra $A_\Delta$ is the quotient of the free
exterior algebra on symbols $e_{\alpha}$ ($\alpha\in\Delta$) by its
ideal generated by the elements
$$
\sum_{j=1}^s (-1)^{j-1} 
e_{\alpha_1}\wedge\cdots\wedge\widehat{e_{\alpha_j}}
\wedge\cdots\wedge e_{\alpha_s}
$$
where $\alpha_1,\ldots,\alpha_s\in\Delta$ are linearly dependent (see
\cite{OT} 3.5). It follows that the space $A_\Delta[r]$ is generated
by the elements
$$
\omega_{(\alpha_1,\ldots,\alpha_r)}:=
\omega_{\alpha_1}\wedge\cdots\wedge\omega_{\alpha_r}=
\frac{d\alpha_1\wedge\cdots \wedge d\alpha_r}{\alpha_1\cdots\alpha_r}
$$
where $(\alpha_1,\ldots,\alpha_r)$ is an ordered basis of $\Delta$.
Moreover, the linear relations between the
$\omega_{(\alpha_1,\ldots,\alpha_r)}$ are consequences of the relations
$$
\sum_{1\leq j\leq r,c_j\neq 0} (-1)^{j-1}
\omega_{\alpha_1,\ldots,\widehat{\alpha_j},\ldots,\alpha_r,\alpha}=0
$$
where $(\alpha_1,\ldots,\alpha_r)$ is as above, and where
$\alpha=\sum_{j=1}^r\,c_j\alpha_j$.

Finally, a basis of $A_\Delta[r]$ consisting of certain
$\omega_{(\alpha_1,\ldots,\alpha_r)}$ can be defined as follows (see
\cite{GZ}, \cite{OT} 3.2 and \cite{S}). Choose an ordering
$(\alpha_1,\alpha_2,\alpha_3,\ldots,\alpha_N)$ of $\Delta$. Consider
the subset $B\subset \CB(\Delta)$ consisting of the ordered bases
$b=(\alpha_{i_1}, \alpha_{i_2},\dots, \alpha_{i_n})$
(listed with strictly increasing indices) such that for all
$j\neq i_p$, the set $\{\alpha_j\}\cup\{\alpha_{i_p};i_p>j\}$
is linearly independent. Then the $\omega_b$ ($b\in B)$ are the
desired basis.

Translating these results in terms of $S_\Delta$ leads to the
following

\begin{proposition}\label{basis}
The set $(\phi_b)_{b\in B}$
is a basis of $S_{\Delta}$. Furthermore, the space 
of linear relations between the $\phi_\sigma$ 
($\sigma\in{\cal B}(\Delta)$) is  generated by the ``Orlik-Solomon
relations''
$$
r_{\sigma,\alpha}:=\phi_{\sigma}-\sum_{\beta\in\sigma}
c_{\alpha\beta}\phi_{\sigma\cup\{\alpha\}\setminus\{\beta\}}
$$ 
where $\sigma\in{\cal B}(\Delta)$, $\alpha\in\Delta\setminus\sigma$
and $\alpha=\sum_{\beta\in\sigma}c_{\alpha\beta}\beta$.
\end{proposition}

For completeness, we will present an a priori proof of this result;
first, let us give an example.

\noindent
{\bf Example.} Let $V$ be a vector space with basis $e_1, e_2, e_3$.
Consider the ordered set
$$
\Delta=(e_1, e_2, e_3,  e_1+e_2, e_2+e_3, e_1+e_2+e_3)
$$
(the set of positive roots of a root system of type $A_3$). Then the
set $B$ consists of  
$$b_1=(e_1,e_2, e_3),~b_2=(e_1,e_3,e_1+e_2),~b_3=(e_1,e_2, e_2+e_3),$$
$$b_4=(e_1, e_1+e_2, e_2+e_3),~b_5=(e_1,e_2,e_1+e_2+e_3),~
b_6=(e_1,e_3,e_1+e_2+e_3).$$
\begin{proof} Let $L$ be the free vector space with basis the elements
$\phi_\sigma$, $\sigma\in{\cal B}(\Delta)$. Let $LR$ be the kernel of
the natural map from $L$ to $R_\Delta$. By definition, $LR$ is the
space of linear relations between the $\phi_\sigma$.
We denote by $C$ the subspace of $L$ with basis 
$(\phi_b)_{b\in B}$, and by $OSR$ the subspace of $LR$
generated by the elements $r_{\sigma,\alpha}$. Let us show that
$L=C+OSR$. If $\sigma=(\alpha_{i_1},\alpha_{i_2},\ldots,\alpha_{i_r})$
is an ordered basis of $\Delta$, we set 
$|\sigma|=i_1+i_2+\cdots+i_r$.
If $\sigma$ is not in $B$,
then there exists a $j$ such that the set
$\{\alpha_j\}\cup\{\alpha_{i_p}; i_p>j\}$ is
linearly dependent. Using the relation 
$r_{\sigma,\alpha_j}$, we replace 
$\phi_{\sigma}$ by a linear combination of  elements
$\phi_{\tau}$ where $\tau$ is 
obtained from $\sigma$ by replacing
one of the elements $\alpha_{i_p}$ with $i_p>j$
by $\alpha_j$. It follows that the numbers $|\tau|$
are strictly smaller than  $|\sigma|$ , so that 
by induction, we obtain $L=C+OSR$.
This shows that the set $(\phi_b)_{b\in B}$
generates $S_\Delta$, and that $LR$ is spanned by $C\cap LR$ and
by $OSR$.

We now show that $C\cap LR=0$.  We need to check that 
if $\sum_{b\in B} c_b \phi_b=0$ as a rational function, then 
all $c_b$ are equal to $0$. We prove this by induction on the number
of elements in $\Delta$.
Remark that all elements of 
$B$ contain $\alpha_1$.
For a set $\kappa$ of $r-1$ linearly independent vectors, 
let $H(\kappa)\subset V$
be the hyperplane generated by $\kappa$. For a
hyperplane $H\subset V$, set 
$\CB(H):=\{\kappa\subset \Delta; H(\kappa)=H\}$.
We write 
$$
\sum_{b\in B} c_b \phi_b=\alpha_1^{-1}\sum_H \phi_H$$
with 
$$\phi_H:=\sum_{b,H(b\setminus\{\alpha_1\})=H}
\frac{c_{b}}{\prod_{\beta\in b\setminus\{\alpha_1\}}\beta}.
$$
Choose a hyperplane $H_0$ generated 
by $b\setminus\{\alpha_1\}$, for some $b\in B$.
Then the residue operator $Res_{V/H_0}$ 
kills all elements $\alpha_1^{-1}\phi_H$
except $\alpha_1^{-1}\phi_{H_0}$, which is mapped to $\phi_{H_0}$. 
Thus, $\phi_{H_0}=0$.
But remark that if we consider the ordered set
$\Delta_0=\Delta\cap H_0$,
the set $B_0=B(\Delta_0)$ consists exactly of the elements $b_0$ 
such that $\{\alpha_1\}\cup b_0\in B$.
We conclude by applying the induction
hypothesis to the vector space $H_0$ and the system $\Delta\cap H_0$,
for all $H_0$.
\end{proof}\bigskip

We see that giving an ordering of $\Delta$, the set
$B$ is characterized as the unique generating family 
$(\phi_b)$ with $\sum_{b\in B}|b|$
minimum.

Following \cite{S}, we now construct the dual basis
$(\phi^b)_{b\in B}$ of the basis $(\phi_b)$, by using iterated
residues. In our framework, they can be introduced as follows.

Let $\alpha\in\Delta$ and let $\Delta \setminus k\alpha$ be the
complement in $\Delta$ of the set of scalar multiples of
$\alpha$. Denote by $\Delta/\alpha$ the image of 
$\Delta\setminus k\alpha$ in the quotient space $V/k\alpha$. Any
$\phi\in S_{\Delta}$ has at worst a simple pole along $\alpha=0$. 
Thus, restriction of $\alpha\phi$ to $(\alpha=0)=(V/k\alpha)^*$ is
a well defined element of $(\Delta/\alpha)^{-1}S(V/k\alpha)$; we
denote it by $Res_{\alpha}(\phi)$. This defines a linear map
$$
Res_{\alpha}:S_{\Delta}=G_{\Delta}[-r]\to 
G_{\Delta/\alpha}[-r+1]=S_{\Delta/\alpha}.
$$
Given an ordered basis $(\beta_1,\ldots,\beta_r)$ of $\Delta$, we
can iterate this construction to obtain a linear form
$$Res_{\beta_1}Res_{\beta_2}\cdots Res_{\beta_r}$$
on $S_{\Delta}$. On the other hand, we have a complete flag
$$0=V_0\subset V_1\subset \cdots\subset V_r=V$$
(where each $V_j$ is spanned by 
$\beta_{r-j+1},\beta_{r-j+2},\ldots,\beta_r$)
together with a non-zero element of each $\wedge^jV_j$. Thus, we
obtain another linear form
$$Res_{V_1/V_0}\cdots Res_{V_{r-1}/V_{r-2}} Res_{V/V_{r-1}}$$
on $S_{\Delta}$ identified with 
$S_\Delta\otimes\wedge^r V$.

\begin{proposition}\label{iterated}
1) (Orlik-Solomon-Terao) For any $\alpha\in\Delta$, the map
$Res_{\alpha}$ defines an exact sequence 
$$
0\to S_{\Delta\setminus k\alpha}\to S_{\Delta}\to
S_{\Delta/\alpha}\to 0.
$$

2) For any ordered basis $b=(\beta_1,\ldots,\beta_r)$ of $\Delta$,
we have
$$Res_{\beta_1}Res_{\beta_2}\cdots Res_{\beta_r}=
Res_{V_1/V_0}\cdots Res_{V_{r-1}/V_{r-2}} Res_{V/V_{r-1}}
$$
with notation as above.

3) (Szenes) The dual basis of $(\phi_b)_{b\in B}$ is given by
$$
\phi^b=Res_{\beta_1}Res_{\beta_2}\cdots Res_{\beta_r}.
$$
\end{proposition}
\begin{proof}
1) is a consequence of \cite{OT} Theorem 3.126; a direct proof is a
follows. Let $\sigma$ be a basis of $\Delta$ containing $\alpha$, and
let $\sigma/\alpha$ be its image in $V/k\alpha$. Then $\sigma/\alpha$
is a basis of $\Delta/\alpha$ and all bases of $\Delta/\alpha$ are
obtained in this way. Moreover, $Res_{\alpha}(\phi_{\sigma})$ is a
non-zero multiple of $\phi_{\sigma/\alpha}$. It follows that
$Res_{\alpha}$ is surjective. 

Clearly, the kernel of $Res_{\alpha}$ contains 
$S_{\Delta\setminus k\alpha}$. Conversely, if 
$\phi\in S_{\Delta}$ 
is mapped to 0 by $Res_{\alpha}$, then $\phi$ is defined on
$(\alpha=0)$ and thus, $\phi\in R_{\Delta\setminus k\alpha}$. We can
write $\phi=\phi_1+\phi_2$ where 
$\phi_1\in S_{\Delta\setminus k\alpha}$ and 
$\phi_2\in V^*R_{\Delta\setminus k\alpha}$. Then
$\phi_2=\phi-\phi_1$ is in $S_{\Delta}$, too, whence $\phi_2=0$
and $\phi\in S_{\Delta\setminus k\alpha}$.

2) Let $b'=(\beta'_1,\ldots,\beta'_r)$ be another ordered basis of
$\Delta$. Consider the element 
$$\rho:=Res_{\beta_1}\cdots Res_{\beta_r}(\phi_{b'}).$$
If $Res_{\beta_r}(\phi_{b'})$ is non-zero, then we must have
$\beta_r=t_r\beta'_{\pi(r)}$
for some non-zero $t_r\in k$ and some index $\pi(r)$. If moreover
$Res_{\beta_{r-1}}Res_{\beta_r}(\phi_{b'})$ is non-zero, then we
must have 
$\beta_{r-1}\in t_{r-1}\beta'_{\pi(r-1)}+k\beta'_{\pi(r)}$
for some non-zero $t_{r-1}\in k$ and some $\pi(r-1)\neq\pi(r)$
(because $\beta_{r-1}$ is not a multiple of $\beta_r$). Continuing
in this way, we see that either $\rho=0$ or there exists a permutation
$\pi$ of $\{1,\ldots,r\}$ and non-zero $t_1,\ldots,t_r\in k$ such that
$$
\beta_i\in t_i\beta'_{\pi(i)}+\sum_{j=i+1}^rk\beta'_{\pi(j)}
$$
for all $i$. Then we have $\rho=t_1\cdots t_r$. 

On the other hand, set
$$
\rho':=Res_{V_1/V_0}\cdots Res_{V/V_{r-1}}(\phi_{b'}).
$$
If $Res_{V/V_{r-1}}(\phi_{b'})\neq 0$, then
there exist a unique index $\pi'(1)$ and a unique non-zero $t'_1\in k$
such that
$\beta'_{\pi'(1)}\in t'_1\beta_1+\sum_{j=2}^r k\beta_j.$
Further, $\beta'_i\in\sum_{j=2}^r k\beta_j$ for all $i\neq\pi'(1)$. 
If moreover 
$Res_{V_{r-1}/V_{r-2}}Res_{V/V_{r-1}}(\phi_{b'})\neq 0$, then 
$\beta'_{\pi'(r-1)}\in t'_2\beta_2+\sum_{j=3}^r k\beta_j$
for uniquely defined $\pi'(r-1)$ and $t'_{r-1}$. Further,
$\beta'_i\in\sum_{j=3}^r k\beta_j$ for all
$i\notin\{\pi'(1),\pi'(2)\}$. Continuing, we obtain if $\rho'\neq 0$:
$$
\beta'_{\pi'(i)}\in t'_i\beta_i+\sum_{j=i+1}^r k\beta_i
$$
for a permutation $\pi'$ and non-zero $t'_1,\ldots,t'_r$; then we have
$\rho'=1/t'_1\cdots t'_r$. This is equivalent to the set of conditions
of the first part of the proof, with $\pi'=\pi$ and $t'_i=1/t_i$.

3) Let $b'=(\beta'_1,\ldots,\beta'_r)\in B$ such that 
$Res_{\beta_1}\cdots Res_{\beta_r}(\phi_{b'})\neq 0$.
Let $\pi$ and $t_1,\ldots,t_r$ be as above; then $\beta_{r-1}$,
$\beta'_{\pi(r-1)}$ and $\beta_r=t_r\beta'_{\pi(r)}$ are linearly
dependent. Write $\beta_{r-1}=\alpha_i$, $\beta'_{\pi(r-1)}=\alpha_{i'}$
and $\beta_r=\alpha_j$, then $j>i$ and $j>i'$. If $i>i'$ then
$\alpha_{i'}$, $\beta_{r-1}$ and $\beta_r$ are linearly dependent,
which contradicts the hypothesis $b\in B$. Similarly, we cannot have
$i<i'$. Thus, $i=i'$, that is, $\beta_{r-1}=\beta'_{\pi(r-1)}$. In
this way we obtain $\beta_q=\beta'_{\pi(q)}$ for all $q$; because $b$
and $b'$ are in $B$, it follows that $b=b'$.
\end{proof}

\section{Laplace transform and Jeffrey-Kirwan residue}

Starting from now, we assume that $k=\R$.
The Laplace transform associates to a polynomial function 
supported on an acute cone in 
$V$ a rational function on 
$V^*$. We will define the inverse 
Laplace transform, formally denoted  by 
$\int_{V^*}^{\delta}e^{\langle y,h\rangle }[[\phi(y)]] dy$,  
of a meromorphic function  $\phi$ on $V^*$ with poles on a set of
hyperplanes. It depends of a choice of a chamber
$\delta$ in $V^*$ and is a locally polynomial function
on the dual cone $\delta^{\vee}$ of $\delta$. 

Let $V$ be an oriented real vector space of dimension $r$.
We denote by $o$ its orientation.
We consider  as before a finite subset 
$\Delta$ of $V\setminus\{0\}$, which spans $V$; we assume moreover
that $-\alpha\in\Delta$ for all $\alpha\in\Delta$.
A wall of $V$ is defined to be an hyperplane generated by 
$r-1$ linearly independent elements of $\Delta$.
We denote by 
$\CH(\Delta)$ the union of walls. This is a set of hyperplanes in $V$.
The set
$$
V_{reg,\Delta}:=V-\CH(\Delta)
$$
is the set of regular elements in $V$.

We define the vector space 
$\CP_{\Delta}$ of locally polynomial functions on 
$V_{reg,\Delta}$. Elements of 
$\CP_{\Delta} $ are given by polynomial functions
on each connected component of  
$V_{reg,\Delta}$.
The space $S(V^*)$ identifies to the space of polynomial functions on $V$.
Thus the space $\CP_{\Delta} $ is a module under the action 
of $S(V^*)$ by multiplication.

Let $C\subset V$ be an acute convex cone with non empty interior.
Let $C^{\vee}\subset V^*$ be its (closed) dual cone. Then the interior
of $C^{\vee}$ is not empty.
We denote by $[C]$ the characteristic function of $C$, that is,
the function with value 1 on $C$ and 0 outside $C$.

Let $f\in S(V^*)$ be a polynomial function on $V$.
Let $dh$ be an element of $\wedge^r V^*$. As $V$ is oriented, we can
integrate over $V$ 
a differential form $\alpha=\phi(h)dh$ of maximal degree
(here $\phi$ is an integrable function on $V$). We denote the integral 
over $V$ of such a differential form $\alpha$ by $\int_{V,o}\alpha$.
A change of orientation produces a change of sign.

For each $y$ in the interior of $C^{\vee}$, the integral
$$L^o(f[C])(y):=\int_{V,o}e^{-\langle y,h\rangle}f(h)[C](h)dh$$
converges, and defines the {\sl Laplace transform} of $f[C]$.
If $C$ is generated by multiples of elements of $\Delta$, it is easy
to see that $L^o(f[C])$ is given by 
the restriction to the interior of $C^{\vee}$ of a rational function
on $V^*$ belonging to the subspace  $G_{\Delta}$ of $R_{\Delta}$.
We still denote this rational function by 
$L^o(f[C])(y)$.
More exactly, as $L^o(f[C])(y)$
depends linearly on $dh$, we see that
$L^o(f[C])(y)$ is a rational function with 
values in $(\wedge^r V^*)^*=\wedge^r V$.
Furthermore it is clear that the map 
$L^o$ interchanges the action of 
$S(V^*)$ by multiplication on $\CP_{\Delta} $
with its action by derivation on $G_{\Delta}$, up to the automorphism
$P(h)\mapsto P(-h)$.

Let $\delta$ be a connected component of 
the set $V^*-\CH^*(\Delta)$.
Then $\delta$ is an open acute polyhedral cone in $V^*$, and 
$\delta^{\vee}$ is a closed acute polyhedral cone in $V$.
We denote 
by ${\cal P}_{\Delta}(\delta)$ the subspace of 
$\CP_{\Delta} $ spanned by functions 
$f(h) [C(\sigma)](h)$
where $f\in S(V^*)$ and 
$\sigma$ is a basis of $\Delta$ {\bf such that}
$C(\sigma)\subset \delta^{\vee}$; here $C(\sigma)$ denotes the closed
convex cone generated by $\sigma$. 

If $dh$ is a positive element of 
$\wedge^{max}V^*$, we denote by 
$vol(\sigma,dh)$ the volume of the parallelelepiped
constructed on the basis $\sigma$ for the positive density $dh$
corresponding to the differential form $dh$. Specifically, if 
$dh=e^1\wedge e^2\wedge\cdots \wedge e^r$
and if $\sigma=\{\alpha_1,\ldots,\alpha_r\}$,
we have $vol(\sigma,dh)=|\det\langle\alpha_i,e^j\rangle_{i,j}|$. 
Finally, we denote by $L^o_{\delta}$
the restriction to ${\cal P}_{\Delta}(\delta)$ of the Laplace
transform $L^o$.
\begin{theorem}\label{iso}
Given any chamber
$\delta$ in $V^*$,
the Laplace transform 
$$
L^o_{\delta}:{\cal P}_{\Delta}(\delta)\to G_{\Delta}\otimes\wedge^{max}V
$$
is an isomorphism and commutes with the actions of $S(V^*)$ up to the
automorphism $P(h)\mapsto P(-h)$.

We have, for $dh\in \wedge^{max}V^*$ positive (with respect to 
our choice of orientation $o$), and $\sigma$ a basis of 
$\Delta$ such that $C(\sigma)\subset \delta^{\vee}$,
$$
<L^o_{\delta}[C(\sigma)],dh>=\vol(\sigma, dh)\phi_{\sigma}
$$
that is, for $y\in \delta$, we have the equality of functions
$$
\int_{C(\sigma),o}e^{-\langle y,h\rangle}dh=\vol(\sigma, dh)\phi_{\sigma}(y).
$$
\end{theorem}
\begin{proof}
The formula for the Laplace transform of $[C(\sigma)]$ is
straightforward. It implies surjectivity of $L^o_{\delta}$ because
this map is $S(V^*)$-linear, and the $S(V^*)$-module $G_{\Delta}$ is
generated by the $\phi_{\sigma}$ where $\sigma$ is a basis of $\Delta$
such that $C(\sigma)$ is contained in $\delta^{\vee}$ (here we use the
assumption that $\Delta$ is centrally symmetric).

For injectivity of $L^o_{\delta}$, we observe that any function 
$\phi\in{\cal P}_{\Delta}(\delta)$ 
is supported in the acute cone $\delta^{\vee}$, and that $\phi=0$ if
and only if $\phi$ vanishes outside a set of measure zero. Moreover,
the set of functions
$h\mapsto e^{-\langle y,h\rangle}$
(where $y\in \delta$) is dense in the space of smooth, rapidly
decreasing functions on $\delta^{\vee}$.
\end{proof}
\bigskip

Consider the inverse $(L^o_{\delta})^{-1}:G_{\Delta}\otimes
\wedge^{max}V\mapsto {\cal P}_{\Delta}(\delta)$.
Via the projection map 
$Princ_{\Delta}$,
we can extend the map 
$(L^o_{\delta})^{-1}$ to $\Delta^{-1}\hat S(V)\otimes \wedge^{max}V$.
Thus we set 
$$
F^o_\delta(\phi\otimes dy):=
(L^o_{\delta})^{-1}(Princ_{\Delta}(\phi)\otimes dy).
$$
Thus, $F^o_\delta$ associates to any 
meromorphic function $\phi$ on 
$V^*$ with poles on the hyperplanes 
$\alpha=0$ a locally polynomial function on 
$\delta^{\vee}$.
We denote $F^o_\delta(\phi\otimes dy)$ by the formal notation:
$$
(F^o_\delta(\phi\otimes dy))(h)=
\int_{V^*,o}^{\delta}e^{\langle y,h\rangle}[[\phi(y)]]dy.
$$

We now show that $F_\delta^o$ commutes with the actions of $S(V)$ by
derivations on ${\cal P}_\Delta(\delta)$, and by multiplication on
$R_\Delta$.

\begin{lemma}\label{der}
For any $\psi\in S(V)$ and $\phi\in G_\Delta$, we have
$$
\psi(\partial)F_\delta^o(\phi\otimes dy)=
F_\delta^o(\psi\phi\otimes dy).
$$
\end{lemma}
\begin{proof}
It is enough to check this for $\psi=v\in V$. Then, for any
$y\in\delta$, we have
$$\displaylines{
\int_V (\partial(v)F_\delta^o(\phi))(h)e^{-\langle y,h\rangle} dh
-\langle v,y\rangle \phi(y)
\hfill\cr\hfill
=\int_V (\partial(v)F_\delta^o(\phi))(h)e^{-\langle y,h\rangle} dh
+\int_V F_\delta^o(\phi)(h)\partial(v)(e^{-\langle y,h\rangle}) dh
\hfill\cr\hfill
=\int_V\partial(v)(F_\delta^o(\phi)(h)e^{-\langle y,h\rangle}) dh
=\int_{\Sigma} 
F_\delta^o(\phi)(h) e^{-\langle y,h\rangle} i_v(dh)
\cr}$$
where $\Sigma$ denotes the boundary of the support of
$F_\delta^o(\phi)$; here the latter equality follows from Stokes'
theorem. Because $\Sigma$ is a union of polyhedral cones of smaller
dimensions, the function 
$$
y\mapsto\int_{\Sigma} 
F_\delta^o(\phi)(h) e^{-\langle y,h\rangle} i_v(dh)
$$
is in $NG_\Delta$. We thus have 
$$
L_\delta^o(\partial(v)F_\delta^o(\phi\otimes dy))-v\phi\in NG_\Delta
$$
which implies our formula.
\end{proof}
\bigskip

If $\phi\in S_{\Delta}$,
the image $F^o_{\delta}(\phi\otimes dy)$
is a locally constant function 
on $V_{reg,\Delta}$.
Thus we obtain a number of residue maps defined by chambers 
$\gamma$ in $V$ and $\delta$ in $V^*$:
$$
Res_{\gamma,\delta}:S_{\Delta}\otimes \wedge^{max}V\to{\bf R},$$
$$
\phi\otimes dy\mapsto(F^o_{\delta}\phi)\vert_{\gamma}.
$$
The formula of Theorem \ref{iso}
determines $Res_{\gamma,\delta}(\phi_\sigma\otimes dy)$
for $C(\sigma)\subset \delta^{\vee}$ and $dy$ a positive element
of $\wedge^rV$.
More precisely, if $dh$ is the dual measure to $dy$,
$$
Res_{\gamma,\delta}(\phi_\sigma\otimes dy)
=\frac{1}{\vol(\sigma, dh)},\hspace{1cm}{\rm if}\,\,\gamma\subset
C(\sigma),
$$
$$
Res_{\gamma,\delta}(\phi_\sigma\otimes dy)
=0,\hspace{1cm}{\rm if}\,\,\gamma\cap  C(\sigma)=\emptyset.
$$

As $F_\delta^o$ commutes with the action
of differential operators with constant coefficients,
we have 
$$
F^o_\delta(P(\partial)\phi_\sigma \otimes dy)
=P(-h)F^o_\delta(\phi_\sigma \otimes dy)(h)
$$
so that if $C(\sigma)\subset \delta^{\vee}$
\begin{equation}\label{g}
F^o_\delta(P(\partial)\phi_\sigma \otimes dy)(h)
=\frac{1}{\vol(\sigma, dh)}P(-h)[C(\sigma)](h).
\end{equation}

\begin{proposition}\label{k}
(Jeffrey-Kirwan) For $\phi\in \hat R_{\Delta}$ and $h\in V$,
we have  
$$
F^o_\delta(\phi \otimes dy)(h)=
F^o_\delta(Res_{\Delta} (e^{h}\phi)\otimes dy).
$$
\end{proposition}
\begin{proof}
It is sufficient to prove this formula 
for $\phi=P(\partial)\phi_{\sigma}$.  
As we have for $y\in V^*$:
$$Res_{\Delta} (e^{h}\partial(y)\phi)
=-Res_{\Delta}((\partial(y)e^h)\phi)
=-\langle y,h\rangle Res_{\Delta}(e^h\phi),$$
we obtain
$$
Res_{\Delta}(e^h P(\partial)\phi_{\sigma})
=P(-h)Res_{\Delta}(e^h\phi_{\sigma})
=P(-h)\phi_{\sigma}.
$$
So we see, from Formula (\ref{g})
above, that  the equation of Proposition \ref{k}
is satisfied.
\end{proof}
\bigskip

Proposition \ref{k} provides an effective tool to compute 
the inverse Laplace transform of a rational function $\phi$ 
with poles on hyperplanes. Indeed,  
the function $Res_{\Delta} (e^{h}\phi)$
is an element of $S_\Delta$ (depending of $h$),
so that it can be written as a linear combination
$$
Res_{\Delta}(e^h\phi)=\sum_{\sigma} c_\sigma(h)\phi_\sigma.
$$
The choice of a chamber $\delta$ determines a sign 
$\epsilon(\sigma,\delta)$
for which  
$$
\phi_\sigma=\epsilon(\sigma,\delta)\phi_{\sigma^{\delta}}
$$
where the cone
$\sigma^{\delta}$  has the same axes as $C(\sigma)$ and is contained
in $\delta^{\vee}$. Thus, the restriction to a chamber $\gamma$ in $V$
of the inverse Laplace transform $F^o_\delta(\phi\otimes dy)$
is obtained by summing the polynomial terms
$\epsilon(\sigma,\delta)c_\sigma(h)\vol(\sigma,dh)^{-1}$
for all $\sigma$ such that $\gamma\subset C(\sigma^{\delta})$:
$$F^o_\delta(\phi\otimes dy)|_\gamma
=\sum_{\sigma,\gamma\subset C(\sigma^{\delta})}
\epsilon(\sigma,\delta) c_\sigma(h)\vol(\sigma,dh)^{-1}.$$
This is Jeffrey-Kirwan algebraic formula.

\noindent
{\bf Example}

Let us consider a two-dimensional vector space $V$
with basis $(e_1,e_2)$.
Let $\Delta=\{e_1, e_2, e_1+e_2\}$. Consider 
$$\phi(z_1,z_2)=
\frac{1}{z_1 z_2 (z_1+z_2)}.$$
We have
$$\displaylines{
Res_{\Delta}(e^{h_1z_1+h_2z_2}\phi(z_1,z_2)\otimes dz_1dz_2)
=\frac{h_1 z_1+ h_2 z_2}{z_1 z_2 (z_1+z_2)}
\hfill\cr\hfill
= \frac{h_1}{z_2 (z_1+z_2)}
+ \frac{h_2}{z_1 (z_1+z_2)}.\cr}
$$

If $\delta$ is the component 
$e_1>0, e_2>0$ of $V^*$, we then obtain the 
following picture for the inverse Laplace transform 
of $\frac{1}{z_1 z_2 (z_1+z_2)}$.


In the next section, we determine the change of $F^o_{\delta}\phi$
when jumping over a wall.


\section{The jump formula}

We consider, as in Section 5, a real oriented vector space 
$(V,o)$ with a system of hyperplanes defined by 
$\Delta\subset V-\{0\}$.
Let $\delta$ be a chamber in $V^*$, and let $F^o_{\delta}$ be 
the inverse Laplace transform.
In this section, we relate the jumps of 
$F^o_\delta(\phi)\otimes dy$
across walls, with the poles of the function $\phi$ along the wall.

Let $(V_0,o_0)$ be an oriented wall with its system 
$\Delta_0=\Delta\cap V_0$.
The wall $V_0$ separates $V$ in two half-spaces. Choose 
an equation $z$ of $V_0$ such that $o=z\wedge o_0$, and define
$$
V_+=\{ h\in V, \langle z,h\rangle>0\},
$$
$$
V_-=\{ h\in V, \langle z,h\rangle<0\}.
$$
If $U$ is a component of $(V_0)_{reg,\Delta_0}$
there exists unique components $U_{\pm}$ of $V_{reg,\Delta}$ 
contained in $V_{\pm}$ and such that
$U\subset \overline{U_{\pm}}$.

Let $f\in \CP_{\Delta}$ be a locally polynomial 
function on $V_{reg,\Delta}$.
Then the restriction of $f$
to $U_+$ (resp. $U_-$)
is given by a polynomial function $f^+$ (resp. $f^-$).
We define the locally polynomial function 
$Jump_{o/o_0}(f)\in \CP_{\Delta_0}$
by the formula
$$Jump_{o/o_0}(f)|_U=f^+|_U-f^-|_U.$$
\begin{theorem}\label{jump}
Let $(V_0,o_0)$ be an oriented wall. Let $\delta$ be a chamber
in $V^*$ and $\delta_0$ a chamber in $V_0^*$
such that $\delta_0^{\vee}\subset \delta^{\vee}$.
Then, for any $\phi\in \hat R_{\Delta}$, we have the Jump formula:
$$
Jump_{o/o_0}(F_\delta^o(\phi\otimes dy_0))=
F_{\delta_0}^{o_0}(Res_{V/V_0}(\phi\otimes dy)).
$$
\end{theorem}
\begin{proof}
It is sufficient to prove this formula 
for $\phi\in G_\Delta$. (On $NG_\Delta$, both sides are equal to $0$,
because $Res_{V/V_0}$ maps $NG_{\Delta}$ to $NG_{\Delta_0}$).
Thus it is sufficient to prove this formula
for a derivative $\phi=P(\partial)\phi_\sigma$ of an element 
$\phi_\sigma$, with $C(\sigma)\subset \delta^{\vee}$. Then 
$$
F_{\delta}^o(\phi\otimes dy)(h)=P(-h)[C(\sigma)](h).
$$

If $V_0$ is not a wall of $C(\sigma)$,
then $F_\delta^o(\phi\otimes dy)$ has no jump along $V_0$.
Thus the left-hand side of the equality in Theorem \ref{jump}
is equal to $0$. The right-hand side is also $0$, as 
there are at least $2$ vectors in $\sigma$
which are not in $\Delta_0$.

If $V_0$ is a wall of $C(\sigma)$, there exists $\beta\in\Delta$
such that $\sigma=\sigma_0\cup\{\beta\}$
where $\sigma_0$ is a basis of $V_0$. Write $V=V_0\oplus \R \beta$.
Write an element $h\in V$ as $h=h_0+h_1 \beta$
with $h_0\in V_0$ and $h_1\in \R$.
Then the left-hand side is the function 
$P(-h_0)[C(\sigma_0)](h_0)$. If $P$ is divisible by $h_1$, then
$Res_{V/V_0}(P(\partial)\phi_{\sigma}\otimes dy)=0$. Thus, both sides
vanish. If $P$ only depends on $h_0$, then
$$
Res_{V/V_0}(P(\partial)\phi_{\sigma}\otimes dy)=
P(\partial)\phi_{\sigma_0}\otimes dy_0
$$ 
whence the right-hand side is $P(-h_0)[C(\sigma_0)](h_0)$.
\end{proof}
\bigskip

As an application of the Jump formula, let us relate the behaviour at
infinity of a function $\phi\in G_\Delta$ to the order of
differentiability of its inverse Laplace transform. For a positive
integer $n$, we say that $\phi$ {\bf vanishes at order $n$ at infinity}
if the rational function $t\mapsto t^{n-1}\phi(y+tz)$ is 0 at $\infty$
for all regular $y\in V^*$ and for all $z\in V^*$. Equivalently,
$\psi\phi\in G_\Delta$ for any $\psi\in S(V)$ of degree at most $n-1$
(indeed, recall that $G_\Delta$ is the space of functions that vanish
at infinity).

\begin{corollary}
For a function $\phi\in G_{\Delta}$ and a non-negative integer $k$,
the following conditions are equivalent:

1) $F_\delta^o(\phi\otimes dy)$ extends to a function of class $C^k$
on $V$.

2) $\phi$ vanishes at order $k+2$ at infinity.

\noindent
Further, for a wall $V_0$ with equation $z_0=0$ and for $\phi$
satisfying (1) or (2), the following conditions are equivalent:

1)' $F_\delta^0(\phi\otimes dy)$ extends to a function of class
$C^{k+1}$ along $V_0$.

2)' For any regular $z\in V^*$, the rational function
$t\mapsto\phi(z+tz_0)$ vanishes at order $k+3$ at infinity.
\end{corollary}
\begin{proof}
Observe that $F_\delta^o(\phi\otimes dy)$ extends to a continuous
function on $V$ if and only if it has no jumps along walls. This
amounts to $Res_{V/V_0}(\phi\otimes dy)=0$ for any wall $V_0$
(because $Res_{V/V_0}$ maps $G_\Delta\otimes\wedge^r V$ to 
$G_{\Delta_0}\otimes \wedge^{r-1}V_0$, and
$F_{\delta_0}^{o_0}$ is injective on the latter). Equivalently,
$$
Res_{t=\infty}(\phi(z+tz_0)dt)=0
$$
for all regular $z$ and for all $z_0$. Because $\phi$ vanishes at
infinity, this means that $\phi$ vanishes at order 2 there. This
proves the equivalence of (1) and (2) in the case where $k=0$.

The general case follows by induction on $k$. Indeed, recall that
$$
\partial(v)F_\delta^o(\phi\otimes dy)=F_\delta^o(v\phi\otimes dy)
=F_\delta^o(Princ_\Delta(v\phi)\otimes dy)
$$
for any $v\in V$. Thus, using the induction hypothesis for $k-1$,
assertion (1) is equivalent to: $\phi$ and $Princ_\Delta(v\phi)$
vanish at order $k+1$ at infinity. Then $v\phi\in G_{\Delta}$ (because
$\phi$ vanishes at order 2 at infinity) and (1) is equivalent to:
$v\phi$ vanishes at order $k+1$ at infinity.

The proof of equivalence of (1)' and (2)' is similar.
\end{proof}

\section{Orlik-Solomon relations and stratified Fourier transform}

We still consider a real vector space $V$ with a finite subset 
$\Delta\subset V\setminus\{0\}$ such that $\Delta$ spans $V$ and
$\Delta=-\Delta$. We fix a Lebesgue measure
$dh$ on $V$ and a chamber $\delta\subset V^*$. Changing slightly
notation, the inverse Laplace transform $F_{\delta}$ associates to any
element of $G_{\Delta}$ a locally polynomial function on
$V_{reg,\Delta}$. In this section, we associate to any element of
$G_{\Delta}$ a {\bf piecewise polynomial} function defined on all of
$V$. This assignement will depend on the choices of a chamber $\delta$
in $V^*$ and of a chamber $\gamma$ in $V$; it will be denoted by
$F_{\gamma,\delta}$. The piecewise polynomial
function $F_{\gamma,\delta}(\phi)$ will extend the locally polynomial
function $F_{\delta}(\phi)$, and will be the continuous extension of
$F_{\delta}(\phi)$ if it exists. We will use the function
$F_{\gamma,\delta}(\phi)$ in part III of this article, in connection
with the definition of Eisenstein series.

Denote by ${\cal PP}_{\Delta}$ the vector space of functions on $V$
spanned by functions $P[C]$ where $P\in S(V^*)$ and $C$ is a closed
polyhedral cone with axes in $\Delta$. Then ${\cal PP}_{\Delta}$ is a
$S(V^*)$-submodule of the module of piecewise polynomial functions on
$V$ (for the stratification where the open strata are the chambers,
and the closures of other strata are proper faces of closures of
chambers). We begin by constructing a morphism of $S(V^*)$-modules
from $G_{\Delta}$ to a quotient of ${\cal PP}_{\Delta}$. This morphism
will depend on the choice of a chamber $\gamma$ in $V$, will be
denoted by $F_{\gamma}$, and will be called the formal Fourier
transform.

For a basis $\sigma$ of $V$, we denote by $|\det(\sigma)|$ the volume
of the parallelepiped constructed on $\sigma$.
We set 
$$
a_\sigma :=|\det(\sigma)|\phi_{\sigma}=\frac{|\det(\sigma)|}
{\prod_{\alpha\in \sigma}\alpha},
$$ 
an element of $S_\Delta$. Remark that $a_\sigma$
does not change if we multiply elements in $\sigma$ by positive
constants. The Orlik-Solomon relations are more naturally expressed in
terms of the $a_\sigma$, as shown by the following result, an easy
consequence of Theorem \ref{free} and Proposition \ref{basis}.
\begin{proposition}\label{OS}
Let $\sigma$ be a basis of $\Delta$, let
$\alpha\in\Delta\setminus\sigma$ and let
$\alpha=\sum_{\beta\in\sigma} c_{\alpha\beta}\beta$ be the expansion
of $\alpha$ in the basis $\sigma$.
Then the elements $(a_\sigma)_{\sigma\in \CB(\Delta)}$
verify the relations
$$a_\sigma=
\sum_{\beta\in \sigma,c_{\alpha\beta}\neq 0}\,sign(c_{\alpha\beta})
a_{\sigma\cup\{\alpha\}\setminus\{\beta\}}.\leqno(OS)$$
Furthermore, if ${\cal M}$ is a $S(V^*)$-module
and $(A_\sigma)_{\sigma\in \CB(\Delta)}$ is a family in
${\cal M}$ verifying the relations (OS), then
there exists a unique map  $A:R_\Delta\to{\cal M}$ such that 

1) the map $A$ commutes with the action of $S(V^*)$.

2) For all $\sigma\in \CB(\Delta)$, we have $A(a_\sigma)=A_\sigma$.

3) $A(NG_\Delta)=0.$
\end{proposition}

Remark that the relations (OS) have coefficients equal to $\pm 1$.
It makes thus sense to find elements in an abelian group,
satisfying these relations.
The group ${\cal C}(V)$ generated by characteristic functions
of locally closed polyhedral cones in $V$ will 
be very useful to construct such elements $A_\sigma$.

We introduce some notation. A {\sl polyhedral cone} in $V$ is a 
closed convex cone $C\subset V$ (with vertex at 0)
which is generated by finitely many vectors $h_1,\ldots,h_n$; we set
$C=C(h_1,\ldots,h_n)$.
For $A$ a subset of $V$, we denote by $[A]$ 
the characteristic function of $A$, i.e.,
the function on $V$ with value 1 on $A$ and $0$
outside $A$. We denote by ${\cal C}(V)$ the additive group of integral
valued functions on $V$, generated by all characteristic functions of
polyhedral cones.

For any closed convex cone $C$, we denote by $C^0$ the relative interior
of $C$, i.e., the interior of $C$ in the affine space generated by $C$.
Observe that ${\cal C}(V)$ contains the characteristic
functions of relative interiors of polyhedral cones, and more generally,
the characteristic functions of locally closed polyhedral cones.
The subgroup of ${\cal C}(V)$ generated by
characteristic functions of polyhedral cones which contain lines is
denoted by ${\cal LC}(V)$. 
For example if $\alpha\in V$ is nonzero, then
$$[C(-\alpha)]+[C(\alpha)^0]\in {\cal LC}(V).$$

We denote by ${\cal C}_{\Delta}$ the subspace of ${\cal C}(V)$
generated by characteristic functions of polyhedral cones $C(\kappa)$ 
where $\kappa\subset\Delta$. Then, by definition, ${\cal PP}_{\Delta}$
is the $S(V^*)$-module generated by ${\cal C}_{\Delta}$. We denote by
${\cal LC}_{\Delta}$ the subspace of ${\cal C}_{\Delta}$
generated by functions $[C(\kappa)]$ where $\kappa\subset\Delta$ and
$C(\kappa)$ contains a line.

Let $p\in V$ and let $C\subset V$ be a polyhedral cone with non-empty
interior, such that $p$ lies in no hyperplane generated by a facet of
$C$. Set
$$
C'_p:=\{h\in C~\vert~{\rm the~segment~} [h,p]{\rm ~meets~} C^0\}.
$$
Then $C'_p$ is equal to $C$ minus the union of its facets which
generate a hyperplane separating $C^0$ and $p$. In particular, $C'_p$
is a locally closed polyhedral cone. If moreover $C=C(\sigma)$ where
$\sigma\in{\cal B}(\Delta)$, and 
$p=\sum_{\alpha\in\sigma}\,p_{\alpha}\alpha$ is in $V_{reg,\Delta}$,
then we obtain easily
$$
C(\sigma)'_p:=
C(\alpha,p_{\alpha}>0)+C(\alpha,p_{\alpha}<0)^0.
$$
In particular, the cone $C(\sigma)'_p$
depends only of the chamber $\gamma$ which contains $p$. Thus, we
denote it by $C(\sigma)'_{\gamma}$.

Define a map
$$
A_{\gamma}:{\cal B}(\Delta)\to {\cal C}_{\Delta}
$$
by
$$
A_{\gamma}(\sigma)=[C(\sigma)'_{\gamma}].
$$
\begin{theorem}
For any chamber $\gamma$ in $V$,
the family of elements $(A_{\gamma}(\sigma))_{\sigma\in \CB(\Delta)}$
verify the relations $(OS)$ in the quotient group 
${\cal C}_{\Delta}/{\cal LC}_{\Delta}$.
\end{theorem}
\begin{proof}
Because of the relation 
$[C(-\alpha)]=-[C(\alpha)^0]$ modulo ${\cal LC}_{\Delta}$
we see that the image of the element $A_{p}(\sigma)$
in ${\cal C}_{\Delta}/{\cal LC}_{\Delta}$
changes sign, if we flip one of the elements $\beta_j$ in 
$\sigma=(\beta_1,\beta_2,\cdots, \beta_r)$ to $-\beta_j$.
We thus may assume that the relation is 
$$\alpha=\beta_1+\beta_2+\cdots +\beta_s$$
for some $s\leq r$. Then the cones 
$C(\sigma\cup\{\alpha\}\setminus\{\beta_j\})$ ($1\leq j\leq s$) are the
maximal cones in a polyhedral subdivision of $C(\sigma)$, and we
conclude by the lemma below.
\end{proof}
\begin{lemma}
Let $C\subset V$ be a polyhedral cone. Let $C_1,\ldots, C_n$ be
the maximal cones of a polyhedral subdivision of $C$. Let $p\in V$
such that $p$ lies in no hyperplane generated by a facet of some
$C_i$. Then $C'_p$ is the disjoint union of 
$C'_{1,p},\ldots,C'_{n,p}$.
\end{lemma}
\begin{proof} Clearly, each $C'_{i,p}$ is contained in $C'_p$.
Conversely, let $x\in C'_p$. If $x$ lies in no $C'_{i,p}$ then the
segment $[x,p]\cap C^0$ has a non-empty interior in $[x,p]$ and is
contained in the union of all facets of the $C_i$. It follows that
this segment is contained in a facet of some $C_i$. Thus, $p$ is in
the hyperplane generated by this facet, a contradiction. So 
$x\in C'_{i,p}$ for some $i$. Assume that $x\in C'_{j,p}$ for some
$j\neq i$. Then $[x,p]\cap C_i^0$ and $[x,p]\cap C_j^0$ are disjoint
segments with non-empty interiors in $[x,p]$. Moreover, because 
$x\in C_i\cap C_j$, the closures of both segments contain $x$, a
contradiction.
\end{proof}
\bigskip

We denote by ${\cal LP}_{\Delta}$ the $S(V^*)$-submodule of 
${\cal PP}_{\Delta}$ generated by ${\cal LC}_{\Delta}$, that is, the
space of piecewise polynomial functions which are polynomial in at
least one direction. By the preceeding theorem, together with Proposition
\ref{OS}, each choice of a chamber $\gamma$ in $V$ defines a morphism
of $S(V^*)$-modules $F_{\gamma}$ from $G_{\Delta}$ to the quotient
space ${\cal PP}_{\Delta}/{\cal LP}_{\Delta}$.
\begin{definition}
Let $\gamma$ be a connected component of $V_{reg,\Delta}$. We denote
by
$$
F_{\gamma}:R_{\Delta}\to{\cal PP}_{\Delta}/{\cal LP}_{\Delta}
$$
the unique map such that

1) $F_{\gamma}$ commutes with the action of $S(V^*)$
up to the automorphism $P(h)\mapsto P(-h)$.

2) $F_{\gamma}(\vert\det(\sigma)\vert\phi_{\sigma})=
[C(\sigma)'_{\gamma}]$
for all $\sigma\in{\cal B}(\Delta)$.

3) $F_{\gamma}(NG_{\Delta})=0$.

\noindent
We call $F_{\gamma}$ the formal Fourier transform.
\end{definition}

Now we construct a lift of 
$F_{\gamma}:R_{\Delta}\to{\cal PP}_{\Delta}/{\cal LP}_{\Delta}$
to ${\cal PP}_{\Delta}$. In other words, we associate to any element
of $R_{\Delta}$ a piecewise polynomial function on $V$, compatibly
with $F_{\gamma}$. We may do this by {\bf specifying a chamber in}
$V^*$, as shown by 

\begin{lemma}
Let $\delta$ be a chamber in $V^*$ and let $\phi\in R_{\Delta}$. 
Then $F_{\gamma}(\phi)$ has a unique representative with support in
$\delta^{\vee}$.
\end{lemma}
\begin{proof}
Let $\alpha\in\Delta$, then $\alpha$ or $-\alpha$ is in
$\delta^{\vee}$. Using the relation 
$[C(-\alpha)]+[C(\alpha)^0]\in{\cal LC}_{\Delta}$, we see that
$C(\kappa)$ has a representative with support in $\delta^{\vee}$, for
any linearly independent $\kappa\subset\Delta$. This shows existence.
For uniqueness, it is enough to check that any 
$f\in {\cal LP}_{\Delta}$ with support in some acute cone $C$ must be
zero. This is shown in the proof of \cite{BV} Theorem 1.4 for 
$f\in{\cal LC}(V)$; this proof adapts with minor changes, as follows.
Embed ${\cal LP}_{\Delta}$ into the vector space ${\cal F}(V^*)$ of
functions on $V^*$. The additive group of $V^*$ acts on 
${\cal F}(V^*)$ by translations; we denote by $z\mapsto T(z)$ this
action. For a polyhedral cone $C$ which contains a line $l$, we have
$(1-T(z))[C]=0$ for all $z\in l$. Thus, for $P\in S(V^*)$, we have
$$
(1-T(z))^N(P[C])=0
$$
whenever $N>deg(P)$. Because $f\in{\cal LP}_{\Delta}$, it follows that
there exist $z_1,\ldots,z_n\in V^*\setminus\{0\}$ (non necessarily
distinct) such that
$$\prod_{j=1}^n(1-T(t_jz_j))f=0$$
for all $t_j\in\R$. Moreover, we can find $h\in V$ such
that $h>0$ on $C\setminus\{0\}$ and that $\langle h,z_j\rangle\neq 0$
for all $j$. Replacing $z_j$ by $-z_j$, we may assume that 
$\langle h,z_j\rangle<0$ for all $j$. Let $w\in V^*$. We can choose
$A>0$ such that 
$$\langle h,w+\sum_{j\in J} t_jz_j\rangle<0$$
for any non-empty subset $J$ of $\{1,\ldots, n\}$ and for $t_j>A$. We
have 
$$0=(\prod_{j=1}^n(1-T(-t_jz_j)f)(w)=
\sum_{J\subset\{1,\ldots,n\}}
(-1)^{\vert J\vert}f(w+\sum_{j\in J}t_jz_j).$$
By assumption, $f$ is identically zero on the open half-space $h<0$.
It follows that $f(w)=0$.
\end{proof}
\bigskip

We denote by $F_{\gamma,\delta}(\phi)$ the representative of
$F_{\gamma}(\phi)$ with support in $\gamma^{\vee}$. Let us compute
$F_{\gamma,\delta}(a_{\sigma})$ for $\sigma\in{\cal B}(\Delta)$.
Write 
$$
a_{\sigma}=\epsilon(\sigma,\delta)a_{\sigma^{\delta}}
$$
where $\epsilon(\sigma,\delta)=\pm 1$ and the cone
$C(\sigma^{\delta})$ has the same axes as $C(\sigma)$ 
and is contained in $\delta^{\vee}$. Then, by definition
$$
F_{\gamma,\delta}(a_{\sigma})=
\epsilon(\sigma,\delta)[C(\sigma^{\delta})]'_{\gamma}.
$$
These elements $F_{\gamma,\delta}(a_{\sigma})$
satisfy the Orlik-Solomon relations in the space ${\cal PP}_{\Delta}$.
\begin{definition}
Let $\gamma$ be a chamber in $V$ and let $\delta$ be a chamber in
$V^*$. We denote by 
$$F_{\gamma,\delta}:R_{\Delta}\to{\cal PP}_{\Delta}$$
the unique map such that

1) $F_{\gamma,\delta}$ commutes with the action of $S(V^*)$
up to the automorphism $P(h)\mapsto P(-h)$.

2) $F_{\gamma,\delta}(\vert det(\sigma)\vert\phi_{\sigma})
=[C(\sigma)'_{\gamma}]$ for all $\sigma\in{\cal B}(\Delta)$ such that
$C(\sigma)\subset\delta^{\vee}$.

3) $F_{\gamma,\delta}(NG_{\Delta})=0$.

\noindent
We call $F_{\gamma,\delta}$ the stratified Fourier transform.
\end{definition}

We now express $F_{\gamma,\delta}(\phi)$ in terms of
$F_{\delta}(\phi)$.
\begin{proposition}
Let $\gamma$ be a chamber of $V$, let $p\in\gamma$ and let $\delta$ be
a chamber in $V^*$. Then we have for any $\phi\in G_{\Delta}$ and
$h\in V$:
$$F_{\gamma,\delta}(\phi)(h)=
\lim_{\epsilon\to 0,\epsilon>0}F_{\delta}(\phi)(h+\epsilon p).$$
In particular, $F_{\gamma,\delta}(\phi)$ is an extension of 
$F_{\delta}(\phi)$ to the whole of $V$, and is the continuous
extension if it exists.
\end{proposition}
\begin{proof}
Observe first that the formula makes sense: because $p$ is regular,
$h+\epsilon p$ is regular for $\epsilon$ sufficiently small and
$\epsilon>0$. If the formula holds for $\phi$ then it holds for
$P(\partial)\phi$ where $P\in S(V^*)$, because both $F_{\delta}$ and
$F_{\gamma,\delta}$ are $S(V^*)$-linear. Thus it suffices to check the
formula for $\phi=\phi_{\sigma}$ where
$C(\sigma)\subset\gamma^{\vee}$. Then 
$F_{\gamma,\delta}(\phi)=[C(\sigma)'_p]$ whereas
$F_{\delta}(\phi)$ is restriction of $[C(\sigma)]$ to
$V_{reg,\Delta}$. But
$$[C(\sigma)'_p](h)=
\lim_{\epsilon\to 0,\epsilon>0}[C(\sigma)](h+\epsilon p)$$
as follows from the definition of $C'_p$.
\end{proof}

\end{document}